\documentclass[12pt]{article}
\usepackage{amsmath}
\usepackage{latexsym}
\usepackage{amssymb}
\usepackage{graphicx}
\newcommand{\mr}{\mathrm{mr}}
\newcommand{\sgn}{\mathrm{sgn}}
\newcommand{\rank}{\mathrm{rank}}
\newenvironment{proof}{\noindent{\bf Proof}\hspace*{1em}}{\qed\bigskip} 
\newtheorem{theorem}{Theorem}

\newcommand{\qed}{\rule{7pt}{7pt}}

\title{\ \\[-0.8in] The Minimum Rank Problem: a counterexample}
\author{ Swastik Kopparty \\
Computer Science and Artificial Intelligence Laboratory \\
MIT \\
Cambridge, MA, U.S.A.\\
E-mail:  swastik@mit.edu\\
and \\
K. P. S. Bhaskara Rao \\
Department of Mathematics and Computer Science \\
Indiana State University\\
Terre Haute, IN 47802 U.S.A. \\
E-mail: bkopparty@isugw.indstate.edu}
\date{}
\begin{document}
\maketitle
\begin{abstract}

We provide a counterexample to a recent  conjecture that the minimum
rank of every sign pattern matrix can be realized by a rational matrix.
We use one of the equivalences of the conjecture and some results from
projective geometry. As a consequence of the counterexample we show that
there is a graph for which the
minimum rank over the reals is strictly smaller than the minimum rank
over the rationals. We also make some comments on the minimum rank of
sign pattern matrices over different subfields of $\mathbb R$.

\end{abstract}

\medskip

\medskip


\section{Introduction}

The main reference for this paper is \cite{AHR} in which the conjecture and its equivalences appear.  

A matrix whose entries are from the set $\{+,-,0\}$ is called a {\it
sign 
pattern matrix}.   A matrix with real entries is called a real matrix and a
matrix with rational entries is called a rational matrix. For a real matrix $B$,
sgn$(B)$ 
is the sign pattern matrix obtained by replacing each positive
(respectively, negative, 
zero) entry of $B$  by $+$ (respectively, $-$, 0). 
If $A$ is a sign pattern matrix and $\mathbb F$ is a subfield of $\mathbb R$, the {\it sign pattern class} of $A$ over
$\mathbb F$  is defined by
 $$Q_{\mathbb F}(A)=\{\ B\ : B \mbox{ is a matrix with entries in $\mathbb F$ and } \mbox{sgn}(B)=A\ \}.$$  

For a sign pattern matrix $A$ and a subfield $\mathbb F$ of $\mathbb R$, the {\it minimum rank} of $A$ over $\mathbb F$, denoted
$\mr_{\mathbb F}(A)$, is defined as 
$$ \mr_{\mathbb F}(A) = \min_{B \in Q_{\mathbb F}(A)} \{ \mbox{rank } B \}.$$

In \cite{AHR}, the authors made the following basic conjecture:
$$ \mbox{ For any $m\times
n$ sign pattern matrix $A$, $\mr_{\mathbb R}(A) = \mr_{\mathbb Q}(A)$.}$$
They showed that the conjecture holds in certain special cases.   
 
In \cite{AHR}, it was also shown that the above conjecture is equivalent to another conjecture, namely, 
\bigskip

\noindent \ \ For any real matrices $D$, $C$, and $E$,
with $DC = E$, there are rational matrices $D^*$, $C^*$, and $E^*$ such
that 
sgn$(D^*) =$ sgn$(D)$, sgn$(C^*) =$ sgn$(C)$, sgn$(E^*) =$ sgn$(E)$, and
$D^*C^* = E^*$.
\bigskip

In the following section, we shall give an example to show that this conjecture is not true.
We also show that there is a graph for which the
minimum rank over the reals is strictly smaller than the minimum rank
over the rationals. In the last section we make some comments on the minimum rank of
sign pattern matrices over different subfields of $\mathbb R$.

\section{The Counterexample}
Consider a configuration $\cal C$ (from \cite{G}, p.92) of nine points and nine lines given by A, B, C, D, E, F, G, H, I, and 
nine lines ABEF, ADG, AHI, BCH, BGI, CEG, CFI, DEI, DFH as drawn in Figure 1 below starting with a regular pentagon. 

Let ${\ell_1, \ell_2, \ldots ,\ell_9}$ be the nine lines in Figure 1 and  let the 
equation of ${\ell_{i}}$ be ${a_{i}x+b_{i}y+c_{i}}=0$. Let the nine points (with real coordinates) 
be $(x_{i}, y_{i})$, $i = 1, 2, \ldots, 9$. 

Let $D$ be the $9 \times 3$ matrix whose $i^{th}$ row is $(a_{i}, b_{i}, c_{i})$  and $C$ be the $ 3 \times 9 $ matrix 
whose $j^{th}$ column is the transpose of the row $(x_{i},y_{i},1)$. 
Let $DC = E$. $E$ is a $9 \times 9$ matrix whose $(i,j)^{th}$ element is $0$ if the $j^{th}$ point is on 
the $i^{th}$ line and   $\neq 0$ if the $j^{th}$ point is not on the $i^{th}$ line. The incidences of the 9 points on the 9 lines 
are exactly dictated by the zero and nonzero elements of $E$. 

The result on p.93 of \cite{G} states that (the incidence structure) $\cal C$ 
cannot be realized with nine points with rational coordinates. Suppose now that there are rational matrices $D^*$, $C^*$, 
and $E^*$ such
that $D^*C^* = E^*$ and the zero non-zero pattern of $E^*$
is same as the zero non-zero pattern of $E$. Since the third row of $C^*$ has nonzero elements, by
dividing each column of $C^*$ and the corresponding column of $E^*$ by a
nonzero rational number
we may assume that the third row of $C^*$ has all $1's$. Now, let
the $j^{th}$ column of $C^*$ be the transpose of
$(x_{j}^*, y_{j}^*, 1)$. If $D^*$ is the $9 \times 3$ matrix whose $i^{th}$
row is $(a_{i}^*, b_{i}^*, c_{i}^*)$, then the $j^{th}$ point
 $(x_{j}^*,y_{j}^*)$ will be on the 
line $a_{i}^* x + b_{i}^* y + c_{i}^* = 0$ if and only if
$(x_{j}, y_{j})$ is on $ \ell_{i}$ for $i = 1, 2, \ldots, 9$. This is
because $D^*C^* = E^*$
and $E^*$ and $E$ have the same zero non-zero pattern.  Hence
$(a_{i}^*, b_{i}^*)$ for $i = 1, 2, \ldots, 9$ will be nine points with rational
coordinates with the same structure  of $\cal C$.

Hence there are no rational matrices $D^*, C^*, E^*$ such that $D^*C^*=E^*$
and $E^*$ has the same zero pattern as $E$. Hence
there are no rational matrices $D^*, C^*$ and $E^*$ such that $D^*C^* = E^*$
and sgn($D^*$) = sgn$(D)$, sgn$(C^*)$ = sgn$(C)$,
sgn$(E^*)$ = sgn$(E)$. 

The above procedure actually gives a real $ 12 \times 12 $ matrix
$B = \left[ \begin{array}{cc}
 I_3 & C \\
 D & E 
 \end{array} \right],$ such that rank$(B)=3$, for which there is
 no rational matrix $F$ such that rank$(F) =3$ and $F$ and $B$ have
the same zero non-zero pattern.
                                                                        
 If 
$A = \left[ \begin{array}{cc}
 0 & B \\
 B^T & 0 
 \end{array} \right],$ then $A$ is a $24 \times 24$ symmetric real matrix for which there is no rational matrix $A^{\ast}$ such that 
sgn$(A)$ = sgn$(A^{\ast})$ and rank$(A)$ = rank$(A^{\ast})$. This in turn gives us a bipartite graph G on $24$ points (with $12$ points on each
side) whose incidence matrix has the zero non-zero pattern of $A$. For this graph the minimum rank over the rationals is strictly more than 
that over the reals (this rank being $6$).

Note that in [1] it was shown that for every
real matrix B of rank 2 there is a rational matrix $F$ of rank 2
such that $B$ and $F$ have the same sign pattern.

\section{General Results}
Incidence structures with properties such as that of Figure 1 were first constructed systematically 
by Maclane~\cite{M} using the ``von Staudt algebra of throws". Theorem 3 of that paper states:

\begin{theorem}[Maclane \cite{M}]
Let $\mathbb K$ be a finite algebraic field over the field of rational numbers. Then there exists a 
matroid $M$ of rank 3 which can be represented by a matrix with elements of $\mathbb K$, while any 
other representation of $M$ by a matrix of elements in a number-field $\mathbb K_1$ requires
$\mathbb K_1 \supset \mathbb K$.
\end{theorem}

Using this theorem along with the argument of the previous section gives us the following general result.
 
\begin{theorem}
Let $\mathbb K$ be a subfield of $\mathbb R$, finite and algebraic over $\mathbb Q$.
Then there exists a sign pattern matrix $A$, such that for any field $\mathbb K_1 \subset \mathbb R$
with $\mathbb K \not\subseteq \mathbb K_1$, $\mr_{\mathbb K}(A) < \mr_{\mathbb K_1}(A)$.
\end{theorem}

In contrast, the situation completely changes for purely transcendental extensions. 

\begin{theorem}
Let $\mathbb F$ be a subfield of $\mathbb R$, and let $\alpha \in \mathbb R$ be transcendental over $\mathbb F$.
Then for any sign pattern matrix $A$, $\mr_{\mathbb F(\alpha)}(A)= \mr_{\mathbb F}(A)$.
\end{theorem}
\begin{proof}
It is clear that for any sign pattern matrix $A$,
$\mr_{\mathbb F(\alpha)}(A) \leq \mr_{\mathbb F}(A)$.
To prove the reverse inequality, it suffices to show that for any
matrix $M$ with entries in $\mathbb F(\alpha)$, there exists a
matrix $M^*$ with entries in $\mathbb F$ such that $\rank(M^*) \leq \rank(M)$
 and $\sgn(M^*) = \sgn(M)$.

Let $M$ be an $m \times n$ matrix with entries in $\mathbb F(\alpha)$.
By multiplying $M$ by a suitable element of $\mathbb F[\alpha]$,
it suffices to prove the theorem when $M$ has entries
in $\mathbb F[\alpha]$ (which is isomorphic to a polynomial ring,
since $\alpha$ is transcendental over $\mathbb F$). For
each $i \in [m],j \in [n]$, let $M_{ij} = P_{ij}(\alpha)$,
where $P_{ij}$ is a polynomial with coefficients in $\mathbb F$.
As $\alpha$ is transcendental, $P_{ij}(\alpha)=0$ if and only
if $P_{ij}$ is the zero polynomial. Thus we may pick
$\beta \in \mathbb F$ sufficiently close to $\alpha$, so that for
each $i,j$, $P_{ij}(\beta)$ has the same sign as $P_{ij}(\alpha)$. Now
let $g: \mathbb F[\alpha] \rightarrow \mathbb F$ be the substitution
homomorphism (of rings) with $g(\alpha) = \beta$.
Define $M^*$ to be the matrix whose $(i,j)$ entry is
$g(M_{ij}) = P_{ij}(\beta)$.

By construction, $\sgn(M^*) = \sgn(M)$. Let $r = \rank(M)$. Consider
any $S \subseteq [m], T \subseteq [n]$ with $|S| = |T| = r+1$. We know
that the $S \times T$ minor of $M$ vanishes. Thus the
determinant $|( \langle M_{ij} \rangle_{i \in S, j \in T})| = 0$. The
corresponding minor of $M^*$ equals the determinant
$|( \langle g(M_{ij}) \rangle_{i \in S, j \in T})|$ which, using the
fact that $g$ is a homomorphism, equals the determinant $g(|(\langle M_{ij} \rangle_{i \in S, j\in T})|)
= g(0) = 0$.
Thus we have shown that any $(r+1)\times (r+1)$ minor of $M^*$ also vanishes, which gives us the result. 
\end{proof}

\begin{figure*}[h]
\centering
\includegraphics[width=1.0\textwidth]{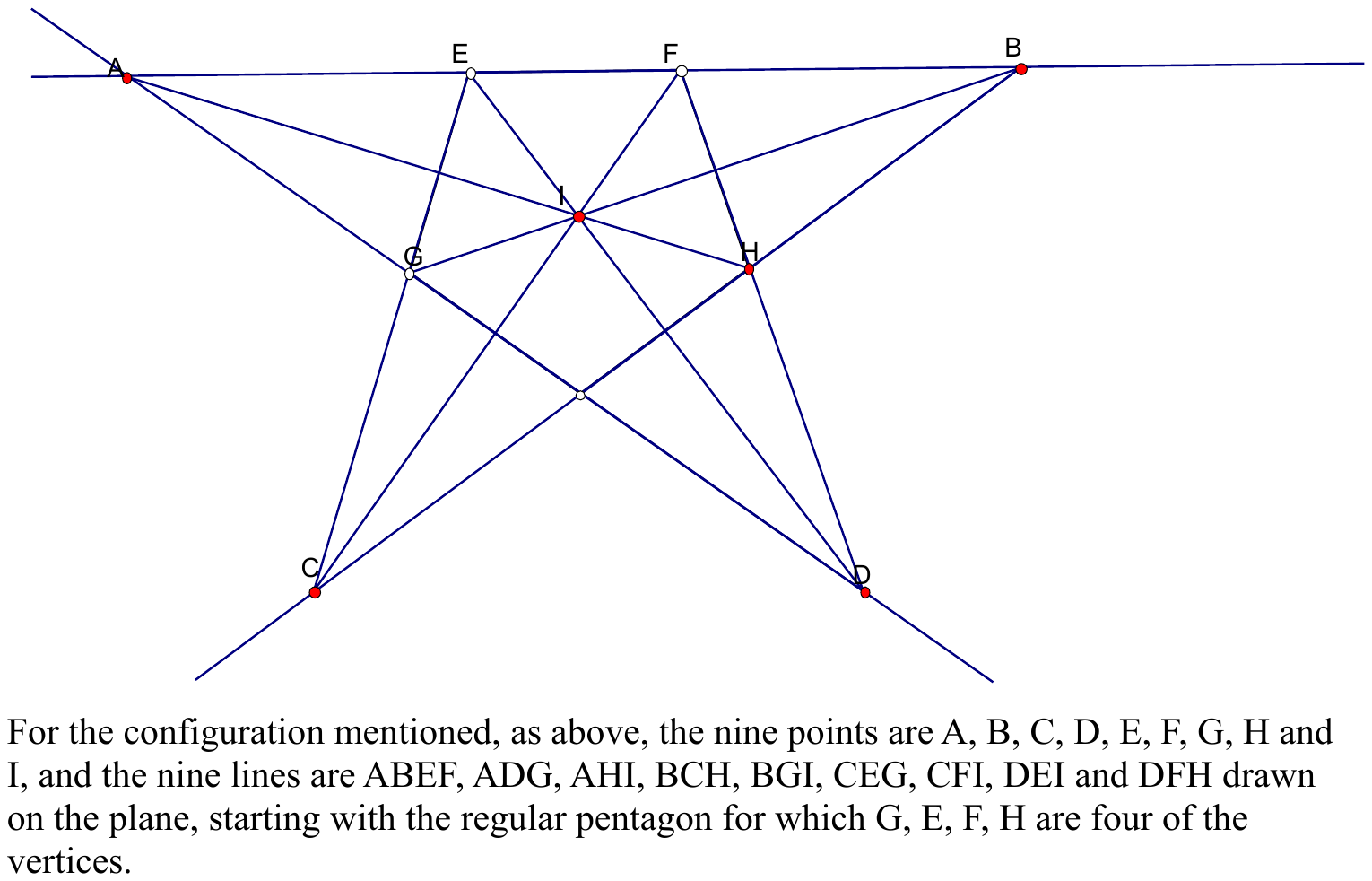}
\caption{}
\label{fig}
\end{figure*}

\end{document}